\documentclass[12pt]{article}
\usepackage{amssymb}
\def \Z{{\mathbb Z}}
\def \p{{\mathbb P}}

\begin{document}
\title{Some remarks on morphisms between Fano threefolds}
\author{Ekaterina Amerik}
\date{April 30, 2004}
\maketitle

Some twenty-five years ago, Iskovskih classified the smooth complex Fano
threefolds with Picard number one. Apart from $\p^3$ and
the quadric, his list includes 5 families of Fano varieties of index two
and 11 families of varieties of index one (for index one threefolds,
 the cube of the anticanonical
divisor takes all even values from 2 to 22, except 20). Recently, the author
(\cite{A}) and C. Schuhmann (\cite{S})
made some efforts to classify the morphisms between such Fano threefolds, the
starting point being a question of Peternell: let $f:X\rightarrow Y$
be a non-trivial morphism between Fano varieties with Picard number one, is
it then true that the index of $X$ does not exceed the index of $Y$?

In particular, Schuhmann (\cite{S}) proved that there are no morphisms from index-two
to index-one threefolds, and that any morphism between index-two threefolds
is an isomorphism (under certain mild additional hypotheses, some of which
were handled later in \cite{A}, \cite{IS}). 
As for morphisms from index-one to index-two Fano  threefolds, such 
morphisms do exist: an index-two threefold has a double covering (branched
along an anticanonical divisor) which is of index one. It is therefore natural
to ask if every morphism from index-one Fano threefold $X$ 
with Picard number one
to index-two Fano threefold $Y$ with Picard number one is a double covering.
In \cite{A}, I proved a theorem (Theorem 3.1) indicating that the answer 
should be yes,
however not settling the question completely. The essential problem was
that the methods of \cite{A} would never work for $Y=V_5$, the linear section
of the Grassmannian $G(1,4)$ in the Pl\"ucker embedding (all smooth 
three-dimensional linear sections of $G(1,4)$ are isomorphic). Though there are
several ways to obtain 
 bounds for the degree of a morphism between
Fano threefolds with second Betti number one (\cite{HM},
\cite{A}), these bounds are still too rough for our purpose.

This paper is an attempt to handle this problem. 
The main result is the following

\

{\bf Theorem} {\it Let $X$ be a smooth complex Fano threefold of index one 
and such that $Pic(X)=\Z$.
Suppose moreover that $X$ is anticanonically embedded. 
Let $f: X\rightarrow V_5$ be a 
non-trivial morphism. Then $X$ is of degree 10 (``$X$ is of type $V_{10}$'') and $f$ is
a double covering. In other words, $X$ is a hyperquadric section of a cone
over $V_5$ in $\p^7$.}

\

I believe that the extra assumption made on $X$ is purely technical and
can be ruled out if one refines the arguments below. This
assumption excludes two families of Fano threefolds: sextic double solids
and double coverings of the quadric branched along a hyperquartic section.
A smooth anticanonically embedded Fano threefold of index one and Picard 
number one is sometimes called {\it a prime Fano threefold}. We shall also
call it thus throughout this paper.

\

{\bf 1. Preliminaries: the geometry of $V_5$}

\

Let us recall some more or less 
classical facts on the threefold $V_5\subset \p^6$, most of which can be found 
in \cite{I} or \cite{FN}.
First of all, as any Fano threefold of index two and Picard number 1, it has
a two-dimensional family of lines. A general line has trivial normal
bundle (call it a $(0,0)$-line), whereas there is a one-dimensional
subfamily of lines with normal bundle 
${\cal O}_{\p^1}(-1)\oplus {\cal O}_{\p^1}(1)$ (call them $(-1,1)$-lines).
The Hilbert scheme of lines on $V_5$ is isomorphic to $\p^2$, the curve
of $(-1,1)$-lines is a conic in this $\p^2$, and there are 3
lines through a general point of $V_5$. More precisely,
the $(-1,1)$-lines form the tangent surface $D$ to a rational normal sextic 
$B$ on 
$V_5$ (in particular, they never intersect), and there are three lines
through any point away from $D$, two lines through a point on $D$ but not on 
$B$, and one line through a point of $B$. The surface $D$ is of degree 10, 
thus 
a hyperquadric section of $V_5$.

We shall denote by $U$ resp. $Q$ the restriction to $V_5$ 
of the universal bundle $U_G$
resp.
the universal quotient bundle $Q_G$ on the Grassmannian $G(1,4)$. 
The cohomology groups related to those bundles
are computed starting from the cohomologies of vector bundles on the 
Grassmannian.
In particular the bundles $U$ and $Q$ remain stable.

We shall also use the following result from \cite{S}:
let $X$ be a prime Fano threelold, and let
$f:X\rightarrow V_5$ be a finite morphism. Let $m$ be such that
$f^*{\cal O}_{V_5}(1)= {\cal O}_X(m)$. Then the inverse image of a general line
consists of $\frac{m^2deg(X)}{10}$ disjoint conics; in general, if one replaces
$V_5$ by another Fano threefold $Y$ of index two with
Picard number one, the inverse image of a general line shall consist 
of $\frac{m^2deg(X)}{2deg(Y)}$ disjoint conics. Here by $deg(Y)$
we mean the self-intersection number of the ample generator of $Pic(Y)$.

Our starting point is the observation that the inverse image of a $(-1,1)$-line
must be connected. This will be the main result of this paragraph.

The Schubert cycles of type $\sigma_{1,1}$ , 
which are sets of points of $G(1,4)$
corresponding to lines lying in a fixed hyperplane, and are also caracterized
as zero-loci of
sections of the bundle dual to the universal, are 4-dimensional
quadrics in the Pl\"ucker $G(1,4)$, so each of them intersects $V_5$ along a 
conic.
Conversely, every smooth conic on $V_5$ is an intersection with such a
Schubert cycle. Indeed, every conics on a Grassmannian
is obviously contained in some $G(1,3)$; and if this conic
is strictly contained in
$G(1,3)\cap V_5$, then  $G(1,3)\cap V_5$ is a surface, so the bundle
$U^*$ has a section vanishing along a surface; but this contradicts the 
stability of $U^*$.
 
The same is (by the same argument) true for pairs of intersecting lines
on $V_5$. Moreover, 
the correspondence between the Schubert cycles and the conics is
one-to-one (it is induced by the restriction map on the global sections 
$H^0(G(1,4),U^*_G)\rightarrow H^0(V_5, U^*)$ which is an isomorphism).

Let us show that among these conics, there is a one-dimensional
family of double lines.

\

{\bf Proposition 1.1} {\it Fix an embedding $V_5\subset G(1,4)\subset \p^9$. 
There
is a one-dimensional family of Schubert cycles $\Sigma_t$ such that 
for each $t$, the intersection of $V_5$ and $\Sigma_t$ is a (double) line.
Moreover, lines on $V_5$ which are obtained as a set-theoretic intersection
with a Schubert cycle of type $\sigma_{1,1}$, are exactly $(-1,1)$-lines.}

\

{\it Proof:} The three-dimensional linear sections of $G(1,4)$ in the Pl\"ucker
embedding are parametrized by the Grassmann variety $G(6,9)$; let, for 
$P\in G(6,9)$, $V_P$ denote the intersection of $G(1,4)$ with the corresponding
linear subspace (which we will denote also by $P$). The Schubert cycles
are parametrized by $G(3,4)=\p^4$; likewise, denote by $\Sigma_t$ the
Schubert cycle corresponding to $t\in \p^4$. Consider the following 
incidence subvariety ${\cal I}\subset G(6,9)\times \p^4$:

$$  {\cal I}= \{(P,t)\in G(6,9)\times \p^4| V_P\cap\Sigma_t \ is\ a\ line\}.$$ 

The fiber ${\cal I}_t$ of ${\cal I}$ over any $t\in \p^4$ parametrizes
the six-dimensional subspaces $P$ of $\p^9$ intersecting $\Sigma_t$ along a 
line. $\Sigma_t$ is a quadric in $\p^5 \subset \p^9$, and $P$ intersects it 
along a line $l$ if and only if the plane $H=P\cap \p^5$ is tangent to
$\Sigma_t$ along $l$, i.e. lies in every $T_x\Sigma_t, x\in l$.
The intersection of all tangent spaces to
$\Sigma_t\subset \p^5$ along $l$ is a three-dimensional projective space 
(the tangent spaces form a pencil of hyperplanes in $\p^5$, because $\Sigma_t$
is a quadric). This means that for every $l$, the planes tangent to
$\Sigma_t$ along $l$ form a one-dimensional family. The family of lines
on a 4-dimensional quadric (= $G(1,3)$) is a 5-dimensional flag variety,
so the planes in $\p^5$ tangent to $\Sigma_t$ along a line are parametrized
by a six-dimensional irreducible variety (a $\p^1$-bundle over a flag
variety). This implies that ${\cal I}_t$ is irreducible of codimension 3
in $G(6,9)$, so ${\cal I}$ is irreducible of codimension 3 in 
$G(6,9)\times \p^4$. 

We must show that the first projection $p_1: {\cal I} \rightarrow G(6,9)$ 
is surjective and its general fiber is of dimension one. First of all, 
remark that
there are points $P$ in the image of $p_1$ such that the corresponding $V_P$
is smooth (so, is a $V_5$). Indeed, fix, as above, $\Sigma_t$, 
$l\subset \Sigma_t$, $H$ a plane in $\p^5=<\Sigma_t>$ such that 
$H\cap \Sigma_t = l$; the remark will
follow if we show that for
a general $\p^6=  P\subset \p^9$ containing $H$, $G(1,4)\cap P$ is smooth. 
We have $H\cap G(1,4)=H\cap \Sigma_t=l$ 
(because $G(1,4)\cap <\Sigma_t>= \Sigma_t$),
so the smoothness
away from $l$ is obvious, and one checks,
again by standart dimension count, that for $x\in l$, the set
$A_x = \{P|H\subset P, G(1,4)\cap P\ is\ singular\ at\ x\}$ is of codimension
two in the space of all $P$'s containing $H$. Therefore for $P$ general in
the image of $p_1$, $V_P$ is smooth.

It is clear that if a smooth $V_P=G(1,4)\cap P$ is such that 
$V_P\cap \Sigma_t=l$, then the corresponding plane $H$ is tangent along $l$
not only to $\Sigma_t$, but also to $V_P$. Thus the normal bundle
$N_{l,V_P}$ has a subbundle $N_{l,H}$ of degree 1, and so $l$ is of type
$(-1,1)$ on $V_P$. Since we have only one-dimensional
family of $(-1,1)$-lines on a smooth $V_P$, we deduce that a fiber of $p_1$ 
over a point $P$ such that
$V_P$ is smooth, is at most one-dimensional. The irreducibility of ${\cal I}$
now implies that $p_1$ is surjective and its general fiber is of dimension
one. This proves the Proposition.

\

Let us now recall the following result of
Debarre (\cite{D}, partial case of Th\'eor\`eme 8.1, Exemple 8.2 (3)):

{\it Let $X$ be an irreducible projective variety, and let 
$f:X\rightarrow G(d,n)$ be a morphism. Let $\Sigma$ be a Schubert cycle
of type $\sigma_m$. If in the cohomologies of $G(d,n)$, 
 $[f(X)]\cdot\sigma_{m+1}\neq 0$, then $f^{-1}(\Sigma)$ is connected.}

Let $X$ be an irreducible projective variety and $f: X\rightarrow V_5$
be a surjective morphism.
Composing with  the embedding $i:V_5\rightarrow G(1,4)$, we can view $f$
as a morphism to $G(1,4)$. The Schubert cycles that we have just
considered are of type $\sigma_{1,1}$ and not $\sigma_2$, however, passing
to the dual projective space, we arrive at the situation of Debarre's theorem
and get the following

\

{\bf Corollary 1.2} {\it For $l$  a $(-1,1)$-line on $V_5$, any irreducible
projective variety $X$ 
and a 
surjective morphism $f:X\rightarrow V_5$,
$f^{-1}(l)$ is connected.}

\

{\bf Remark 1.3} If we knew that the inverse image of a {\it general} line is
always connected, this would immediately solve our problem; indeed, for a Fano
threefold $X$ of index and Picard number one, the equality
$\frac{m^2deg(X)}{10}=1$ implies that $m=1$, $deg(X)=10$ and $f$ is a double
covering. However,
as shows an example of Peternell and Sommese, this is
false in general, even if one supposes that $X$ is a Fano threefold. 
In the example
of \cite{PtS}, $X$ is the universal family of lines on $V_5$, 
which turns out to be
a Fano threefold (of Picard number two, of course), and $f$ is the natural
triple covering. The inverse image of a general line has two connected 
components. 

\

{\bf Remark 1.4} One can ask if there is a similar connectedness
statement for other Fano threefolds of Picard number one and index two. 
Recall that these are
the following: intersection of two quadrics in $\p^5$; cubic in $\p^4$;
double covering of $\p^3$ branched in a quartic; double covering
of the cone over Veronese surface branched in a hypercubic section. 

Smooth quadrics in $\p^5$ are Grassmannians $G(1,3)$, and a smooth intersection
of two quadrics in $\p^5$ is a quadric line complex. 
It is classically known
(see \cite{GH}, Chapter 6) that on a quadric line complex, there is a finite
(and non-zero) number of lines obtained as set-theoretic intersection with
a plane in $G(1,3)$. These lines are obviously $(-1,1)$-lines, since the
corresponding plane is tangent to the quadric line complex along this line.
Our intersection of two quadrics is contained in a pencil of such 
Grassmannians, so there is a one-dimensional family of lines on it such that
each line is the intersection with a plane lying on some Grassmannian of the
pencil. The curve of $(-1,1)$-lines is irreducible (it follows from the
results in \cite{GH}, Chapter 6, that it is smooth and that it is an ample 
divisor on the
Fano surface of lines, in particular, it is connected). Thus it is just the 
closure of that family. So that it follows again from Debarre's paper that
the inverse image of a general $(-1,1)$-line is connected.

As for the cubic, even if such a connectedness statement could hold,
it would not, as far as I see, follow from any well-known general result. 
One can, though, remark that in the examples of Peternell-Sommese type
``(universal family of lines on $Y$)$\rightarrow Y$'', the inverse image of
a $(-1,1)$-line has a tendency to be connected, whereas the inverse image of 
a $(0,0)$-line
is certainly not connected. 
Indeed, it is observed in the literature that, on the
threefolds as above (the cubic, the quadric line complex, $V_5$), a 
line $l$ is in the closure of the curve 
$C_l=$\{{\it lines intersecting
$l$ but different from $l$}\} on the Hilbert scheme if and only if 
$l$ is a $(-1,1)$-line.   

\

{\bf 2. A Hilbert scheme argument}

\

The previous considerations show that on our Fano threefold $X$, a disjoint
union of conics degenerates flatly to a connected l.c.i. scheme. Recall
the following classical example: if one degenerates a disjoint union of two
lines in the projective space into a pair of intersecting lines, the pair of
intersecting lines shall have an embedded point at the intersection. So if
one wants the limit to be a connected l.c.i., this limit must be a double line.
This suggests to ask if a similar phenomenon can occur in our situation, that
is: can it be true that a connected l.c.i. limit of disjoint conics is
necessarily a multiple conic? 

In any case it is easily checked that, say, a connected limit of pairs of 
disjoint conics does not have to have embedded points when the
two conics become reducible and acquire a common component. So even if 
a statement like this could be true, it is probably difficult to prove. In this
paragraph we shall prove a weaker statement: the inverse image of a 
sufficiently general $(-1,1)$-line is either a multiple conic, or supported
on a union of lines.

Let $T$ be the Hilbert scheme of lines on $V_5$ and let ${\cal M}\subset
T\times V_5$ be the
universal family. We have the ``universal family of the inverse images of
lines under $f$'' $${\cal S} = {\cal M}\times_{V_5}X\subset T\times X.$$
Since $f$ is flat and ${\cal M}$ is flat over $T$, ${\cal S}$ is flat over
$T$.

Let $H'$ be the Hilbert scheme of conics on $X$. Consider the irreducible
components of
$H'$ which are relevant for our problem, that is, the components such 
that their 
sufficiently general points correspond to conics
which are in the inverse image of a sufficiently general line on $V_5$.
Denote by $H$ the union of all such components.

For every point $x\in H$, the image of the corresponding conic $C_x$ is a line.
Indeed, ``$f(C)$ is a line''  is a closed condition on conics $C$ because
$f$ is a finite morphism (for $f$ arbitrary, ``$f(C)$ is contained in a line''
would be a closed condition on $C$).

This allows to construct a morphism $p:H\rightarrow T$ taking every
conic to its image under $f$. Indeed, 
$${\cal L}=\{(C, f(x))|x\in C, C\in H\}\subset H\times V_5$$
is a family of lines over $H$; though apriori it is not clear that it is flat,
this is a ``well-defined family of algebraic cycles'' in the sense of Kollar
(\cite{K}, Chapter I)
and so corresponds to a morphism from $H$ to the Chow variety of lines on 
$V_5$, and this is the same as $T$.

We claim that $p$ is finite. Indeed, it is clear that the only obstruction
to the finiteness of $p$ could be the existence of infinitely many double
structures of arithmetic genus zero on some lines on $X$ (``non-finiteness
of the Hilbert-Chow morphism for the family of conics on $X$''). This obviously
happens if one considers conics in $\p^3$ rather than conics on $X$. 
In our
situation, however, this is impossible, and the Hilbert-Chow morphism is even one-to-one.
Indeed, by \cite{I}, the normal bundle of a line in a prime Fano threefold
is either ${\cal O}_{\p^1}\oplus {\cal O}_{\p^1}(-1)$, or 
${\cal O}_{\p^1}(1)\oplus {\cal O}_{\p^1}(-2)$, and there is the following

\

{\bf Lemma 2.1} {\it Let $l\subset X$ be a line on a prime Fano threefold.
If $N_{l,X}={\cal O}_{\p^1}\oplus {\cal O}_{\p^1}(-1)$, then there is no
locally Cohen-Macaulay double structure of arithmetic genus 0 on $l$. If
$N_{l,X}={\cal O}_{\p^1}(1)\oplus {\cal O}_{\p^1}(-2)$, then such a structure
is unique.}

\

{\it Proof:} All locally Cohen-Macaulay double structures on smooth curves
in a threefold are obtained by a construction due to Ferrand (see for example
\cite{BF}, or else \cite{N} for details): if $Y\subset V$ is a smooth curve
on a smooth threefold, and $\tilde{Y}$ is a double structure on $Y$, write
$L$ for ${\cal I}_Y/{\cal I}_{\tilde{Y}}$; in fact $L$ is a locally free 
rang-one ${\cal O}_Y$-module and ${\cal I}_{\tilde{Y}}$ contains 
${\cal I}_Y^2$. The double structure is thus determined by the natural 
surjection
from the conormal bundle of $Y$ in $V$ to $L$, up to a scalar. 
Now take $Y=l$, $V=X$ and let
$L$ be as above; we have an exact sequence
$$0\rightarrow L\rightarrow {\cal O}_{\tilde{l}}\rightarrow {\cal O}_l
\rightarrow 0,$$ from which it is clear that $p_a(\tilde{l})=0$ if and only
if $L={\cal O}_{\p^1}(-1)$. Now in the first part of our assertion, 
there is no 
non-trivial surjection from $N_{l,X}^*$ to ${\cal O}_{\p^1}(-1)$, and in the
second part, such a surjection is unique up to a scalar.

\

Note that we do not have to consider curves which are not locally 
Cohen-Macaulay,
since, for example, the above argument shows that there are no higher genus
locally Cohen-Macaulay double structures, and an embedded point decreases
the genus.

Thus, for any $t\in T$, $p^{-1}(t)$ is a finite set $\{h_1, \dots, h_k\}$,
and to each $h_i$ there corresponds one conic $C_i$ on $X$, mapped to
$l_t$ by $f$. The next step is to show that $f$ and $p$ ``agree with
each other'':

\ 

{\bf Lemma 2.2}
{\it Let $t\in T$ be any point and $l_t\in V_5$ be the corresponding line.
Let $h_1, \dots , h_k$ be the points of $p^{-1}(t)$ and $C_1, \dots, C_k$
the corresponding conics on $X$. Then the support of $f^{-1}(l_t)$ is
$\bigcup_i C_i$}. 

\

{\it Proof:} Indeed, for a general $t\in T$, it is true:
$f^{-1}(l_t)=\bigcup_i C_i$. For a special $t\in T$, choose a curve 
$V\subset T$ through $t$, such that $t$ is the only ``non-general'' 
point of $V$ in the
above sense, and let $U=p^{-1}(V)$. Denote by 
${\cal C}_U\subset U\times X$ the restriction to $U$ of the universal family
of conics over $H$. The support of the fiber over $t$ of 
$(p\times id)({\cal C}_U)\subset V\times X$
is equal to $\bigcup_i C_i$. But the family ${\cal S}|_V$ coincides with
 $(p\times id)({\cal C}_U)$ except at $t$. ${\cal S}|_V$ being
flat, it must be the scheme-theoretic closure of 
$(p\times id)({\cal C}_U)|_{V-\{t\}}$
in $V\times X$, and thus the support of ${\cal S}|_V$ is 
$(p\times id)({\cal C}_U)$,
q.e.d. 

\ 
 
Let now $t\in T$ be a point corresponding to a sufficiently general 
$(-1,1)$-line. We know that $f^{-1}(l_t)$ is connected. Suppose that the number
$k$ from the Lemma is $>1$, so that there are several conics in the 
$Supp(f^{-1}(l_t))$. Decompose the set of those conics 
into two disjoint non-empty subsets $\Sigma_1$ and $\Sigma_2$.

\

{\bf Proposition 2.3} {\it There exists a conic in $\Sigma_1$ which has a
common component with a conic in $\Sigma_2$; in other words,
$(\bigcup_{C\in \Sigma_1}C)\bigcap (\bigcup_{C\in \Sigma_2}C)$ cannot be
zero-dimensional.}

\

{\it Proof} Choose a suitable small 1-dimensional disc $(V,0)$ centered
at $t$. The inverse image $p^{-1}V$ is a disjoint union of two analytic
sets $U_1$ and $U_2$ ($U_i$ consists of points corresponding to conics
near those of $\Sigma_i$). Repeat the procedure of the previous lemma: consider
the universal families ${\cal C}_i$ of conics over $U_i$ and their images
${\cal S}_i=(p\times id)({\cal C}_i)\subset V\times X$. Let ${\cal S}^0$,
${\cal S}_i^0$ denote the restriction of our families ${\cal S}$,
${\cal S}_i$ to the punctured
disc $V^0=V-\{0\}$. The family ${\cal S}^0$ is just the disjoint union of
${\cal S}_i^0$. Now take the closure of all those (as analytic spaces) 
in $V\times X$: the closure
of ${\cal S}^0$ is just ${\cal S}|_V$, by flatness, and the closure ${\cal S}_i'$ 
of
${\cal S}_i^0$ has the same support as ${\cal S}_i$, is contained
in ${\cal S}|_V$ and is flat over $V$. 
The fiber of ${\cal S}_i'$ over $0$, denoted $S_i$,
is contained in the fiber $S$ of ${\cal S}$, 
since the tensor multiplication preserves
the surjectivity. So $f^{-1}(l_t) = S$ contains $S_1
\cup S_2$. By construction, $S_i$ are flat limits of disjoint unions of 
$a_i$ conics and
$S$ is a flat limit of disjoint unions of $a_1+a_2$ $(=\frac{m^2deg(X)}{10})$
conics. 

If $S_1$ and $S_2$ do not have common components, then, since by 
flatness $deg(S)=deg(S_1)+deg(S_2)$, this implies $S=S_1 \cup S_2$,
because $S$ is purely one-dimensional (being an inverse image of a line under
a finite morphism). But then we can apply the exact sequence

$$0\rightarrow {\cal O}_S \rightarrow {\cal O}_{S_1}\oplus {\cal O}_{S_2}
\rightarrow {\cal O}_{S_1\cap S_2}\rightarrow 0$$ and get a contradiction, 
since
by flatness $\chi({\cal O}_S)=\chi({\cal O}_{S_1})+\chi({\cal O}_{S_2})$, 
$S_1\cap S_2$ is non-empty and it is zero-dimensional by assumption.
Thus $S_1$ and $S_2$ must have common components, and, as $S_i$ is supported
on $\bigcup_{C\in \Sigma_i}C$, the Proposition is proved.

\

{\bf Corollary 2.4} {\it In the situation as above, $f^{-1}(l_t)$ is supported
either on a single conic, or on a union of lines.}

\

Indeed, the proposition shows that if $f^{-1}(l_t)$ contains more than one
conic, then any conic from $f^{-1}(l_t)$ must have a common component with
the rest of these conics, that is, it must be singular.

\

Some results from commutative algebra allow to prove a stronger
(``local'') version of Proposition 2.3:

\

{\bf Proposition 2.5} {\it In the situation of Proposition 2.3, through
each intersection point $P$ of $\bigcup_{C\in \Sigma_1}C$ and 
$\bigcup_{C\in \Sigma_2}C$ passes some common component of 
$\bigcup_{C\in \Sigma_1}C$ and 
$\bigcup_{C\in \Sigma_2}C$.}

\

{\it Proof:} The family ${\cal S}$ is flat over $T$ which is smooth, and
the fibers are l.c.i., thus locally Cohen-Macaulay. 
It follows (\cite{EGA}, 6.3.1,
6.3.5) that ${\cal S}$ is locally Cohen-Macaulay, and that the same it true 
for the restriction of ${\cal S}$ to any smooth curve in $T$. Suppose that 
Proposition 2.5 is not true for some intersection point $P$. Let 
$x=(t,P)\in T\times X$ be
the point corresponding to $P$ in ${\cal S}$. Consider the restriction
of ${\cal S}$ to a general curve through $t$, and an analytic neighbourhood
of $x$ in this restriction. Clearly, if one removes $x$, this neighbourhood
becomes disconnected: there are at least two branches corresponding
to $Supp{\cal S}_i$ as in Proposition 2.3. But this is impossible
by Hartshorne's connectedness (\cite{H}), which implies that a connected
Cohen-Macaulay neighbourhood
 remains connected if one removes a subvariety
of codimension at least two.

\

{\bf Remark 2.6} The argument of the Proposition is more or less the following:
``if we have a disjoint union of certain smooth curves $A$ and $B$, 
which degenerates flatly
into a certain connected $C$ in such a way that $A$ and $B$ do not 
acquire common
components in the limit, then $C$ will have embedded points at the 
intersection points of the limits of $A$ and $B$, so this is impossible if
we know that $C$ is purely one-dimensional''. Examples show that one cannot 
say 
anything
reasonable
if one allows $A$ and $B$ to acquire common components. But in fact 
our ``$C$'',
that is, $f^{-1}(l_t)$, is more than just purely one-dimensional: it is a 
locally complete intersection. I do not know if its being a flat limit
of  disjoint unions of conics can impose stronger restrictions on its geometry.

\

To illustrate how we shall apply this, let us handle the case
when $f^{-1}(l_t)$ is supported on a single conic. 

\

{\bf Proposition 2.7} {\it In this case $X=V_{10}$ and $f$ is a 
double covering.}

\

{\it Proof:} As the degree of the 
subscheme
$f^{-1}(l_t)$ of $X$ is $\frac{m^2deg(X)}{5}$, this conic is of multiplicity
$\frac{m^2deg(X)}{10}$ in $f^{-1}(l_t)$. That is, the local degree of $f$ near
a general point of such a conic is also $\frac{m^2deg(X)}{10}$. Now this is
the local degree of $f$ along a certain divisor, because we have chosen
the line $l_t$ to be ``sufficiently general among the $(-1,1)$ lines'': it 
varies in a one-dimensional family. This divisor is thus a component of the
ramification divisor of $f$, and $\frac{m^2deg(X)}{10}-1$ is its ramification
multiplicity.

Now the ramification divisor of $f$ is an element of $|{\cal O}_X(2m-1)|$,
and so the local degree of $f$ at its general point is at most $2m$, and if 
it is $2m$, then the ramification divisor is the inverse image of the
surface covered by the $(-1,1)$-lines and set-theoretically a hyperplane
section of $X$. So we have:

$$\frac{m^2deg(X)}{10}\leq 2m,\ mdeg(X)\leq 20,$$
and if the equality holds, then $f$ is unramified outside the inverse image of
the surface of $(-1,1)$-lines. Also,
$\frac{m^2deg(X)}{10}$ must be an integer. The inequality thus only holds
for $deg(X)=10$ and $m=1$ (this is a double covering) or $m=2$ (in this case
it is an equality), and for $deg(X)=4$ and $m=5$ (also an equality).
Let us exclude the last two cases. If $f$ is unramified 
outside the inverse image of the surface of $(-1,1)$-lines, then $p$ is
$\frac{m^2deg(X)}{10}$-to-one everywhere except over the conic parametrizing
the $(-1,1)$-lines on $T=\p^2$. It is thus a topological covering of the
complement to this conic in $T$. But the latter is simply-connected; so that
$H$ has $\frac{m^2deg(X)}{10}$ irreducible components and each one maps 
one-to-one on $T$. Notice that the number $\frac{m^2deg(X)}{10}$ is superiour
to three in both cases. But this is impossible. Indeed, on $V_5$ one has only
3 lines through a general point; whereas, if $H$ has $k$ components, each
component would give at least one conic through a general point of $X$.
Those conics are mapped to {\it different} lines through $f(x)$, because they
intersect; thus $k\leq 3$.

\

{\bf 3. Proof of the Theorem}

\

We have seen that the inverse image of a general $(-1,1)$-line
is supported either on one conic, or on a union of lines, and settled the
first case in the end of the second section. Let us now settle the remaining
case, using Proposition 2.5. 

The following lemma is standart (and follows e.g. from the 
arguments of \cite{M}, Chapter 3):

\

{\bf Lemma 3.1} {\it Let $g: X_1\rightarrow X_2$ be a proper morphism of 
complex quasiprojective varieties, which is finite of degree $d$. Suppose
that $X_2$ is smooth. Then the inverse image of any point $x\in X_2$
consists of $d$ points at most, and if there are exactly $d$ points in the 
inverse
image of all $x\in X_2$, then topologically $g$ is a covering.}

\

Let $H$ be as in the last section, and let ${\cal C}$ be the universal family
of conics over $H$. Each conic of $H$ is contained in the inverse
image of some line on $V_5$, and set-theoretically such an inverse image
is a union of conics of $H$. Denote by $D$ the surface covered by
$(-1,1)$-lines on $V_5$. Recall that through each point of the complement
to $D$ in $V_5$ there are three lines, that $D$ is a tangent surface to a
rational normal sextic and that there are two lines (one 
$(-1,1)$-line and one $(0,0)$-line) through any point of $D$ away from this
sextic and a single line through each point of the sextic.
Since the inverse image of a general $(0,0)$-line is a disjoint union 
of conics of $H$, there are three conics of $H$ through a general point of
$X$, and at least three through any point away from $f^{-1}(D)$.
The natural morphism $q: {\cal C}\rightarrow X$ is proper and finite of
degree three. By the Lemma, there are exactly three conics of $H$ through any 
point
of $X$ away from $f^{-1}(D)$.

Let $l$ be a general $(-1,1)$-line on $V_5$. Consider the case when
 $Z=f^{-1}(l)$ is a 
set-theoretic union of degenerate conics $C_1, ..., C_k$ of $H$.

\

{\bf Lemma 3.2} {\it $Z$ contains a line which belongs to a single $C_i$
(say $C_1$).}

\

{\it Proof:} Suppose the contrary, that is, that any component of $Z$ is 
contained in at least two conics of $H$. Through a general point $x$ of
this component there is at least one more conic of $H$, coming from the
inverse image of the $(0,0)$-line through $f(x)$. This implies that
the morphism $q: {\cal C}\rightarrow V_5$ is three-to-one outside
an algebraic subset $A$ of codimension at least two in $X$. That is, 
${\cal C}-q^{-1}(A)$ is, topologically, a covering of $X-A$. But $X-A$ is
simply-connected because $X$ is Fano and thus simply-connected. This means
that ${\cal C}$ is reducible, consists of three components and each of them
maps one-to-one to $X$. Since $X$ is smooth, it must be isomorphic to 
each of those components (by Zariski's Main Theorem). But this is impossible
because the components are fibered in conics and $X$ has cyclic Picard group.

\

Before continuing our argument, let us recall some well-known facts
on lines on prime Fano threefolds  (\cite{I}). Lines on our Fano threefold 
$X$ are parametrized by
a curve, which may of course be reducible or non-reduced. Its being reduced
or not influences the geometry of the surface covered by lines on $X$.
Namely, if a component of the Hilbert scheme of lines on $X$ is reduced,
then the natural morphism from the correspondent component of the
universal family to $X$ is an immersion along a general line; and there is
a classical computation (\cite{I}, \cite{T}) which says that if its image $M$
is an element of $|{\cal O}_X(d)|$, then a general line of $M$ intersects
$d+1$ other lines of $M$. If a component of the Hilbert scheme of lines
is non-reduced, then the surface $M$ covered by the corresponding lines
is either a cone (but this can happen only on a quartic), or a tangent
surface to a curve. One knows only one explicit example of a Fano threefold
as above such that the surface covered by lines on it is
 a tangent surface to a curve, it is constructed by Mukai and Umemura 
(having been overlooked by Iskovskih)
and has degree 22. The surface itself is a hyperplane section of this threefold
and its lines never intersect.
 
The following Proposition, due to Iliev and Schuhmann,
is the main result of \cite{IS}  slightly reformulated: 

\

{\bf Proposition 3.3} {\it Let $X$ be a prime Fano threefold, ${\cal L}$ a
complete one-dimensional family of lines on $X$ and $M$
the surface on $X$ covered by lines of ${\cal L}$. If $X$ is different 
from the
Mukai-Umemura threefold, then a general line of ${\cal L}$  
intersects at least one
other line of ${\cal L}$.} 

\

{\it An outline of the proof:} If not, then, by what we have said above, 
the surface $M$ must 
be a tangent surface
to a curve. Studying its singularities, Iliev and Schuhmann prove that it
must be a hyperplane section of $X$. Then they show, by case-by-case analysis
(of which certain cases appear already in \cite{A}),
that the only prime
Fano threefold containing a tangent
surface to a curve as a hyperplane section, is the Mukai-Umemura threefold. 

\

``Lines contained in a single $C_i$'' cover a divisor on $X$ as
$Z$ varies 
(this is the branch divisor of $q$).
Since $(-1,1)$-lines on $V_5$ never intersect,
Proposition 3.3 implies that if $X$ is not the Mukai-Umemura threefold, 
then in $Z$ there are at least two lines contained in a single
conic (say, $l_1\subset C_1$ and $l_2\subset C_2$), and that they intersect,
say at the point $P$. Notice that $C_1$ is necessarily different from $C_2$:
otherwise we get a contradiction with Proposition 2.5 by considering
$\Sigma_1=\{l_1\cup l_2\}$, $\Sigma_2$ the set of all the other $C_i$ and the
intersection point $P$.

\

{\bf Claim 3.4} {Both $C_1$ and $C_2$ are  pairs of lines intersecting at 
the point $P$, and $Z$ is supported on $C_1\cup C_2$.  Thus $Z$ is, 
set-theoretically, the union of three or four lines
through $P$.}

\

{\it Proof:}

1) If $C_1$ is a double line, we get a contradiction with Proposition 2.5 by 
considering $\Sigma_1=\{C_1\}$ and the point $P$; the same is true for $C_2$.

2) Let $C_1=l_1\cup l_1'$. If $l_1'$ does not pass through $P$, we get the
contradiction in the same way, thus $P\in l_1'$. Also, $P\in l_2'$, where
$C_2=l_2\cup l_2'$.

3)There are two possibilities: 

a)If $l_1'\neq l_2'$, then there must be another conic from $Z$ through $P$, 
 containing $l_1'$. Indeed, otherwise we again get a contradiction
with Proposition 2.5. In the same way, there is a conic from $Z$ through $P$ 
which
contains $l_2'$. In fact it is the same conic, because
otherwise there are at least four conics through $P$, contradicting Lemma 3.1.
Denote it by $C_3$. No other conic from $Z$ passes through $P$. So
$C_3=l_1'\cup l_2'$, and $l_1', l_2'$ are not contained in conics others
than $C_1, C_2, C_3$.

b) If $l_1'= l_2'$, then no other conic from $Z$ contains this line
(otherwise through its general point there will pass at least four conics
from $H$, the fourth one coming from the inverse image of the correspondent
$(0,0)$-line).

4) Now the union $C_1\cup C_2\cup C_3$ in the case a), resp. the union
$C_1\cup C_2$ in the case b), cannot have any points in common
with the other components of $Z$; otherwise, taking 
$\Sigma_1=\{C_1, C_2, C_3\}$, resp. $\Sigma_1=\{C_1, C_2\}$, 
we obtain a contradiction with Proposition 2.5.
But $Z$ is connected, so $Z$ is supported on the lines $l_1, l_1', l_2,
l_2'$, q.e.d.

\

We are now ready to finish the proof of the theorem stated in the introduction.

\

{\it Proof of the theorem:}
If $X$ is the Mukai-Umemura threefold, then the lines on $X$
never intersect at all, so that $f^{-1}(l)$ must be supported on a single
conic. Proposition 2.7 shows that a morphism from $X$ to $V_5$ is impossible.
(It should be, however, said at this point that the paper \cite{HM} 
contains a better
proof of the non-existence of morphisms from the Mukai-Umemura
threefold onto any other smooth variety, besides $\p^3$!).

If $X$ is not the Mukai-Umemura threefold and $f^{-1}(l)$ is not supported
on a single conic, then we know by Claim 3.4 how $f^{-1}(l)$ looks. Remark
that $f^{-1}(D)$ is a reducible divisor: it has two components,
one swept out by the lines $l_1$ and $l_2$ as $Z$ varies, another by 
$l_1'$ and $l_2'$. Neither component is a hyperplane section: indeed,
if a hyperplane section of $X$ is covered by lines, then it is either
a cone (impossible in our situation), or a general line intersects
two other lines on the surface by the classical
computation from \cite{T} mentioned above, since a hyperplane section cannot 
be a 
tangent surface to a curve by \cite{IS}. Let $k$ be the multiplicity 
of the component corresponding to $l_i$ and $k'$ be the multiplicity of the
component corresponding to $l_i'$. As $f^*(D)$ is a divisor from 
$|{\cal O}_X(2m)|$, $k+k'\leq m$. At the same time, $Z$ must be of degree
$\frac{m^2deg(X)}{5}$, and thus $2k+2k'=\frac{m^2deg(X)}{5}$, so 
$m^2deg(X)\leq 10$, leavng the only possibility $m=1, deg(X)=10$.

\

{\bf 4. Concluding remarks}

\

In this section, we shall make a further (minor) precision on Theorem 3.1 
from \cite{A}.

In that theorem, it was proved that if 
$X$, $Y$ are Fano threefolds with Picard number one and very ample generator 
of the Picard group, $X$ is of index one, $Y$ is of index two different from
$V_5$ (that is, $Y$ is a cubic or a quadric line complex), 
and $f: X\rightarrow Y$ is a surjective morphism, then $f$ is a ``projection'',
that is, $f^*{\cal O}_Y(1)={\cal O}_X(1)$. The argument of the theorem
also worked for $Y$ a quartic double solid, whereas there were some problems
(hopefully technical ones) for $Y$ a double Veronese cone and for $X$ not 
anticanonically embedded.

Even in the ``good'' cases, the theorem proves a little bit less than one
would like; that is, we want $f$ to be a double covering and we prove only
that $f^*{\cal O}_Y(1)={\cal O}_X(1)$. This still leaves the following 
additional possibilities:

(1) If $Y$ is a cubic, $X$ can be $V_{12}$, $deg(f)=4$ ($X$ cannot be $V_{18}$
because of the Betti numbers: $b_3(V_{18})<b_3(Y)$);

(2) If $Y$ is an intersection of two quadrics, $X$ can be $V_{16}$, $deg(f)=4$
(here $V_{12}$ is impossible since in this case the inverse image of a general 
line would consist of $3/2$ conics).

The first possibility can be excluded by using an inequality of \cite{ARV}: 
it says that for a finite morphism $f:X\rightarrow Y$ and a line bundle $L$
on $Y$ such that $\Omega_Y(L)$ is globally generated, 
$deg(f)c_{top}\Omega_Y(L)\leq c_{top}\Omega_X(f^*L)$, so, for $X$ and $Y$ of
dimension three, $deg(f)(c_3(\Omega_Y)+c_2(\Omega_Y)L+c_1(\Omega_Y)L^2)$
must not exceed $c_3(\Omega_X)+c_2(\Omega_X)f^*L+c_1(\Omega_X)f^*L^2$.

Consider the situation
of (1): we may take $L={\cal O}_Y(2)$, and we know that $c_3(\Omega_Y)=6$
and $c_3(\Omega_X)=10$. Using the equalities $c_2(X)c_1(X)=c_2(Y)c_1(Y)=24$,
we arrive at 
$4(6+24-24)\leq 10+48-48$, which is false. So the case (1) cannot occur.

This inequality does not work in the case (2): indeed, now $c_3(\Omega_Y)=0$,
$c_3(\Omega_Y)=2$ and the inequality reads as follows: 
$4(0+24-32)\leq 2+48-64$, so does not give a contradiction. However we can rule
out this case by our connectedness argument. Indeed, the inverse image
of a general $(-1,1)$-line is connected (Remark 1.4) and the inverse image
of a general $(0,0)$-line consists of two disjoint conics. The results
of Section 2 apply, of course, to our situation; it follows that the
inverse image of a general $(-1,1)$-line is either a double conic, or a union
of two reducible conics which have a common component. In both cases, it
is clear that the ramification locus of $f$ projects {\it onto} the surface 
covered
by $(-1,1)$-lines. But the ramification divisor is a hyperplane section
of $V_{16}$, and thus can project onto a surface from $|{\cal O}_Y(4)|$ at 
most. Whereas it is well-known (and follows for example from the
results in \cite{GH}, Chapter 6) that the surface covered by 
$(-1,1)$-lines on $Y$ is an element of  $|{\cal O}_Y(8)|$.

All this put together gives the following

\

{\bf Theorem 4.1} {\it Let $X$, $Y$ be smooth complex Fano threefolds of 
Picard number one,
$X$ of index one, $Y$ of index two. Assume further that the ample generators
of $Pic(X)$ and $Pic(Y)$ are very ample. Then any morphism from $X$ to $Y$
is a double covering}.

\

I would like to mention that the verification of this statement without
the very ampleness hypothesis amounts to a very small number of particular
cases; for instance, if $Y$ is a double Veronese cone, then already the
formula of \cite{ARV} combined with the knowledge of Betti numbers implies
that for any morphism $f:X\rightarrow Y$ with $X$ Fano of index one with cyclic
Picard group,
$deg(f)=2$ and $X$ is a sextic double solid. It seems that one could be able 
to work out the remaining cases without any essentially new ideas.


\begin{thebibliography}{}

\bibitem[A]{A} E. Amerik, Maps onto certain Fano threefolds,  Doc. Math.  2  
(1997), 195--211.

\bibitem[ARV]{ARV} E. Amerik, M. Rovinsky, A. Van de Ven, A boundedness
theorem for morphisms between threefolds, Ann. Inst. Fourier 49 (1999),
405-415.

\bibitem[BF]{BF} C. B\u anic\u a, O. Forster, 
Multiplicity structures on space curves, The Lefschetz centennial conference, 
Part I (Mexico City, 1984),
Contemp. Math., 58 (1986), 47-64. 



\bibitem[D]{D} O. Debarre, Th\'eor\`emes de connexit\'e pour les produits 
d'espaces projectifs et les grassmanniennes,  Amer. J. Math.  118  (1996),  
no. 6, 1347--1367.

\bibitem[EGA]{EGA} A. Grothendieck, Elements de la g\'eom\'etrie alg\'ebrique,
IV, seconde partie, Publ. Math. IHES 24 (1965).

\bibitem[FN]{FN} M. Furushima, N. Nakayama: The family of lines on the Fano 
threefold $V_5$.  Nagoya Math. J.  116  (1989), 111--122.

\bibitem[GH]{GH} Ph. Griffiths, J. Harris, Principles of Algebraic
Geometry, Wiley, 1978.

\bibitem[H]{H} R. Hartshorne, Complete intersections and connectedness,
Amer. J. Math. 84 (1962), 497-508.

\bibitem[HM]{HM} J.M. Hwang, N. Mok, Finite morphisms onto Fano manifolds of
 Picard number 1 which have rational curves with trivial normal bundles,  
J. Algebraic Geom.  12  (2003),  no. 4, 627--651.

\bibitem[I]{I} V.A. Iskovskih, Fano threefolds I, II, Math. USSR Izv. 11
(1977), 485--527, and 12 (1978), 469-506.

\bibitem[IS]{IS} A. Iliev, C. Schuhmann: Tangent scrolls in prime Fano
threefolds, Kodai Math. J. 23 (2000), no.3, 411-431.

\bibitem[K]{K} J. Kollar, Rational curves on algebraic varieties,
Springer-Verlag, 1996.

\bibitem[M]{M} D. Mumford, Algebraic geometry I. Complex projective 
varieties. Springer-Verlag, 1976.

\bibitem[N]{N} S. Nollet, The Hilbert scheme of degree three curves,
Ann. Sci. ENS 30 (1997), 367-384.

\bibitem[PS]{PtS} Th. Peternell, A.J. Sommese, Ample vector bundles and
branched coverings II, preprint math.AG/0211220.

\bibitem[S]{S} C. Schuhmann, Morphisms between Fano threefolds,  
J. Algebraic Geom.  8  (1999),  no. 2, 221--244.

\bibitem[T]{T} A.N. Tyurin, Five lectures on three-dimensional varieties, 
Uspekhi Mat. Nauk  27  (1972), no. 5 (167), 3--50; translation in 
Russian Math. Surveys  27 (1972), no. 5, 1--53. 

\end{thebibliography}
\end{document}